\documentclass{amsart}
\usepackage{amssymb}
\usepackage{graphicx}
\usepackage{amscd}

\newtheorem{theorem}{Theorem}
\theoremstyle{plain}

\newtheorem{corollary}{Corollary}

\newtheorem{definition}{Definition}

\newtheorem{lemma}{Lemma}

\newtheorem{proposition}{Proposition}
\newtheorem{remark}{Remark}

\numberwithin{equation}{section}

\input{tcilatex}

\begin{document}
\title[wavelets and polyanalytic functions]{Super-wavelets versus
poly-Bergman spaces}
\author{Luis Daniel Abreu}
\address{CMUC, Department of Mathematics of University of Coimbra, School of
Science and Technology (FCTUC) 3001-454 Coimbra, Portugal }
\email{daniel@mat.uc.pt}
\urladdr{http://www.mat.uc.pt/\symbol{126}daniel/}
\thanks{This research was partially supported by CMUC/FCT and FCT
post-doctoral grant SFRH/BPD/26078/2005, POCI 2010 and FSE.}
\date{September, 2009}
\thanks{}

\begin{abstract}
Motivated by potential applications in multiplexing and by recent results on
Gabor analysis with Hermite windows due to Gr\"{o}chenig and Lyubarskii, we
investigate vector-valued wavelet transforms and vector-valued wavelet
frames, which constitute special cases of super-wavelets, with a particular
attention to the case when the analyzing wavelet vector is constituted by
functions $\Phi _{n}$ such that $\mathcal{F}\Phi _{n}(t)=t^{\frac{1}{2}}%
\emph{l}_{n}^{0}(2t)$, where $\emph{l}_{n}^{0}$ is a Laguerre function. We
construct an isometric isomorphism between $L^{2}(%
\mathbb{R}
^{+},\mathbf{C}^{n})$ and poly-Bergman spaces, with a view to relate the
sampling sequences in the poly-Bergman spaces to the wavelet frames and
super-frames with the windows $\Phi _{n}$. One of the applications of the
theory is a proof that $b\ln a<2\pi (n+1)$ is a necessary condition for the
(scalar) wavelet frame associated to the $\Phi _{n}$ to exist. This is the
first known result of this type outside the setting of analytic functions
(the case $n=0$, which has been completely studied by Seip in 1993).
\end{abstract}

\maketitle

\section{Introduction}

\subsection{Motivation}

Over the recent years, an increasing attention has been given to the study
of vector-valued versions of time-frequency and time-scale methods, driven
in part by potential applications in \emph{multiplexing}, an important
method in telecommunications, computer networks and digital video.

Roughly speaking, "multiplexing" means: encoding $n$ independent signals $%
f_{k}$ as a single signal $\mathbf{f}$ which contains all the information of
all of the signals $f_{k}$.

Of course, the signals can always be combined into a single one by simple
superposition; what makes the problem nontrivial, from a mathematical point
of view, is to assure that the superposition is done in such a way that,
when required, each of the signals $f_{k}$ can be recovered from the
multiplexed signal $\mathbf{f}$ (\emph{demultiplexing})$\mathbf{.}$ A recent
mathematical approach to this problem has been proposed by means of the
theory of \emph{frames}. A sequence of functions $\{e_{j}\}_{j\in I}$ is
said to be a frame in a Hilbert space $H$ if there exist constants $A$ and $%
B $ such that, for every $f\in H$,
\begin{equation}
A\left\Vert f\right\Vert _{H}^{2}\leq \sum_{j\in I}\left\vert \left\langle
f,e_{j}\right\rangle \right\vert ^{2}\leq B\left\Vert f\right\Vert _{H}^{2}.
\label{frame}
\end{equation}%
In \cite{Balan}, Balan has given an interpretation of standard multiplexing
methods from the point of view of \emph{super-frames}, which are
vector-valued versions of frames (where the Hilbert space $H$ is $L^{2}(%
\mathbb{R}
,%
\mathbb{C}
^{n})$). The concept of super-frames has been investigated from the
perspective of functional analysis by Han and Larson in \cite{DL} and
several results have been obtained concerning the possibility of extending a
frame to a super-frame \cite{Gu}, something that, not surprisingly, usually
requires oversampling by a rate equal to dimension $n$, which represents the
number of multiplexed signals \cite{Oversampling}, \cite{CharlyYurasuper}.

Within super-frames, super-wavelets have received a particular attention in
\cite{DL}, \cite{Gu} and \cite{Oversampling}. Alternative approaches to the
construction of vector-valued wavelets have also been considered recently
\cite{MRA}, \cite{vvACHA}. In this paper we will provide yet another
approach to the construction of vector-valued wavelets with emphasis in
important special cases where, perhaps in a surprising manner, a new
connection with complex analysis is revealed.

To better give a context to our work, let us review two fundamental ideas
which have arisen recently in different mathematical communities.

The first fundamental idea comes from Gabor and complex analysis: if one
considers the super-frames built from shifting and modulating a window with
n Hermite functions, we are lead to a very structured situation, where an
intriguing connection to complex analysis has been explored by Gr\"{o}chenig
and Lyubarskii in \cite{CharlyYurasuper} and then in \cite{Abreusampling},
where it is also shown that the \emph{vector valued} Gabor representations
with vectorial Hermite windows correspond to \emph{polyanalytic} functions
(we provide definitions some lines below), in much the same way the scalar
Gabor representations with Gaussian correspond to analytic functions.

The second fundamental idea is implicit in the link between the results of
Hutnik \cite{Hutnikcr}, \cite{Hutnikm} and of Vasilevski \cite{VasiBergman},
given the intriguing ressenblance between the description of \emph{wavelet
spaces} and \emph{poly-Bergman spaces}.

In the present paper, we describe a connection between vector-valued
wavelets and polyanalytic functions, in analogy to the one existing for the
Gabor case \cite{Abreusampling}, \cite{Abreustructure}. Indeed, polyanalytic
Bergman spaces are the function-theoretic analogues of wavelet spaces, when
one chooses as analysing wavelet a certain vector defined in terms of
Laguerre functions. We will use this setting to provide a theoretical
solution of the problem of multiplexing and demultiplexing vectorial signals
in $L^{2}(%
\mathbb{R}
^{+},%
\mathbb{C}
^{n})$. We remind that, by the Paley-Wiener theorem, this space is unitarily
equivalent to the vectorial Hardy space on the upper-half plane, $H^{2}(%
\mathbf{U},%
\mathbb{C}
^{n}).$

\subsection{Description of the results}

To describe our results we need some definitions.

Let $\mathbf{U=}\{z\in
\mathbb{C}
:\func{Re}z\geq 0\}$ stand for the upper half plane. The \emph{polyanalytic
Bergman space}, $\mathbf{A}^{n}(\mathbf{U})$ (poly-Bergman space, for
short), is constituted by the complex valued functions defined on the upper
half plane and such that
\begin{equation*}
\left( \frac{d}{d\overline{z}}\right) ^{n}f(z)=0\text{ and\ }\int_{\mathbf{U}%
}\left\vert f(z)\right\vert ^{2}dz<\infty \text{.}
\end{equation*}%
The space $\mathbf{A}^{0}(\mathbf{U})=A(\mathbf{U})$ is the usual (analytic)
Bergman space in the upper half plane.

The potential applications of polyanalytic Bergman spaces in multiplexing
are suggested by the orthogonal decomposition of the space $\mathbf{A}^{n}(%
\mathbf{U})$\ in subspaces called \emph{true polyanalytic Bergman spaces},
which will be denoted by $\mathcal{A}^{k}(\mathbf{U})$:%
\begin{equation*}
\mathbf{A}^{n}(\mathbf{U})=\mathcal{A}^{0}(\mathbf{U})\oplus ...\oplus
\mathcal{A}^{n-1}(\mathbf{U}).
\end{equation*}%
In face of this decomposition, if we start with a vector $\mathbf{f=(}%
f_{0},...,f_{n-1}\mathbf{)\in }L^{2}(%
\mathbb{R}
^{+},\mathbf{C}^{n})$ and find a unitary transform (call it $\emph{Ber}^{n}$%
) such that%
\begin{equation*}
f_{k}\in L^{2}(%
\mathbb{R}
^{+})\overset{\emph{Ber}^{k}}{\rightarrow }F_{k}\in \mathcal{A}^{k}(\mathbf{U%
})\text{,}
\end{equation*}%
then we can store all the $F_{k}%
{\acute{}}%
s$ into a single multiplexed function%
\begin{equation*}
F=F_{0}+...+F_{n-1}\in \mathbf{A}^{n}(\mathbf{U}).
\end{equation*}%
Then, thanks to the orthogonal decomposition of $\mathbf{A}^{n}(\mathbf{U})$%
, it is possible to recover each of the $F_{k}$, projecting the multiplexed
signal $F$ into $\mathcal{A}^{k}(\mathbf{U})$. This is possible to be done
using the reproducing kernel $K^{k}(z,w)$ of the spaces $\mathcal{A}^{k}(%
\mathbf{U})$: Given $F\in \mathbf{A}^{n}(\mathbf{U})$, its true polyanalytic
component $F_{k}\in \mathcal{A}^{k}(\mathbf{U})$ can be recovered by the
orthogonal projection of $F$ over the space $\mathcal{A}^{k}(\mathbf{U})$:%
\begin{equation*}
F_{k}(z)=\left\langle F(w),K^{k}(z,w)\right\rangle _{\mathbf{A}^{n}(\mathbf{U%
})}.
\end{equation*}%
One of the main goals of the present paper is the construction of the
unitary transform $\emph{Ber}^{n}:L^{2}(%
\mathbb{R}
^{+})\rightarrow \mathcal{A}^{n-1}(\mathbf{U})$ mentioned above. We will
call it the \emph{true polyanalytic Bergman transform} and define it in
terms of the (analytic) Bergman transform
\begin{equation*}
Ber\text{ }f(z)=\int_{0}^{\infty }t^{\frac{1}{2}}f(t)e^{izt}dt\text{,}
\end{equation*}%
which is unitary $Ber:L^{2}(%
\mathbb{R}
^{+})\rightarrow A(\mathbf{U})$. Starting from this, we introduce the \emph{%
true polyanalytic Bergman transform}
\begin{equation*}
\emph{Ber}^{n}f(z)=\frac{1}{(2i)^{n}n!}\left( \frac{d}{dz}\right) ^{n}\left[
s^{n}Ber\text{ }f(z)\right] \text{.}
\end{equation*}%
It is at this point that the connections to wavelets shows up. Actually, the
key in proving the unitarity of $\emph{Ber}^{n}$ is to recognize that,
defining the functions $\Phi _{n}$ as
\begin{equation*}
\mathcal{F}\Phi _{n}(t)=t^{\frac{1}{2}}\emph{l}_{n}^{0}(2t)\text{,}
\end{equation*}%
where $\emph{l}_{n}^{0}$ is a Laguerre function (see section 2 for
definition), the true polyanalytic Bergman transform is related to the
wavelet transform by the formula
\begin{equation*}
\emph{Ber}^{n}f(z)=s^{-1}W_{\Phi _{n}}(\mathcal{F}^{-1}f)(x,s)\text{.}
\end{equation*}%
The next goal is to define a unitary isomorphism $\mathbf{Be}r^{n}:L^{2}(%
\mathbb{R}
^{+},\mathbf{C}^{n})\rightarrow \mathbf{A}^{n}(\mathbf{U})$. For $\mathbf{f}%
\in L^{2}(%
\mathbb{R}
^{+},\mathbf{C}^{n})$, the \emph{polyanalytic Bergman transform }is
\begin{equation*}
\mathbf{Be}r^{n}\mathbf{f}=\sum_{0\leq k\leq n-1}\emph{Ber}^{k}f_{k}\text{.}
\end{equation*}%
This provides a unitary operator $\mathbf{Be}r^{n}:L^{2}(%
\mathbb{R}
^{+},\mathbf{C}^{n})\rightarrow \mathbf{A}^{n}(\mathbf{U})$, which relates
vector-valued functions to polyanalytic functions. Here vector-valued
wavelet transforms enter the picture. We will interpretate the polyanalytic
Bergman transform is a special case of a general vector valued wavelet
transform defined as%
\begin{equation*}
\mathbf{W}_{\mathbf{g}}\mathbf{f}(x,s)=\sum_{0\leq k\leq
n-1}W_{g_{k}}f_{k}(x,s)\text{,}
\end{equation*}%
where $W_{g_{k}}f_{k}$ stands for the scalar wavelet transform of the kth
component of the vector $\mathbf{f=(}f_{0},...,f_{n-1}\mathbf{)}$, analyzed
with the wavelet $g_{k}$, which is a component of the vector-valued
analyzing wavelet $\mathbf{g}$. The components of the \ vector $\mathbf{g=(}%
g_{0},...,g_{n-1}\mathbf{)}$ are selected in such a way their Fourier
transforms satisfy
\begin{equation*}
\left\langle \mathcal{F}g_{i},\mathcal{F}g_{j}\right\rangle _{L^{2}(%
\mathbb{R}
^{+},t^{-1})}=\delta _{i,j}\text{.}
\end{equation*}%
In particular, one can choose the vector-valued analyzing wavelet $\mathbf{%
\Phi }_{n}=(\Phi _{0},...,\Phi _{n-1})$. In doing so, we arrive at the
relation between the vector-valued wavelet transform and the polyanalytic
Bergman transform:
\begin{equation*}
\mathbf{W}_{\mathbf{\Phi }_{n}}(\mathcal{F}^{-1}\mathbf{f)}(x,s)=s\mathbf{Ber%
}^{n}\mathbf{f}\text{.}
\end{equation*}

Then, we compute explicitly the reproducing kernels of the spaces $\mathcal{A%
}^{k}(\mathbf{U})$, required for the orthogonal projections. They are given
by the formula
\begin{equation*}
K^{n}(z,w)=\frac{1}{(2i)^{n}n!}\left( \frac{d}{dz}\right) ^{n}\left[
s^{n}\Omega _{n}\left( \frac{z+u}{\eta }\right) \right] \text{,}
\end{equation*}%
where%
\begin{equation*}
\Omega _{n}(z)=4(n+z-i)\left( \frac{1}{z+i}\right) ^{3}\left( \frac{z-i}{z+i}%
\right) ^{n-1}.
\end{equation*}

In the final section we study sampling sequences in polyanalytic Bergman
spaces. This is is equivalent to the study of wavelet frames with the
functions $\Phi _{n}.$ The section contains a contribution to the subject of
\emph{affine density}, the study of the density of sequences $\Lambda $
yielding inequalities of the form (wavelet frames):
\begin{equation*}
A\left\Vert f\right\Vert _{H^{2}(\mathbf{U})}^{2}\leq \sum_{(x,s)\in \Lambda
}\left\vert \left\langle f,T_{x}D_{s}\psi \right\rangle \right\vert ^{2}\leq
B\left\Vert f\right\Vert _{H^{2}(\mathbf{U})}^{2}\text{.}
\end{equation*}%
There has been a considerable research activity around the topic of affince
density, and currently there are three different aproaches (see \cite%
{HeilKutyniok}, \cite{HeilKutyniok2}, \cite{KutyniokCA}, \cite{WSun}, \cite%
{Seip2} and the monograph \cite{KutLN}). The definitions of density in Seip
\cite{Seip2} and \cite{HeilKutyniok}, \cite{HeilKutyniok2}, \cite{KutyniokCA}%
, agree in the case of the hyperbolic lattice $\Lambda
(a,b)=\{(a^{m}bk,a^{m})\}$, where the density is $1/b\log a$, up to
constants independent of $a$ and $b$.

It is straightforward to see that there is no universal necessary lower
bound for an arbitrary function to generate a frame associated to $\Lambda
(a,b)$, in contrast to the situation in Gabor analysis \cite{RS}, \cite{Rief}%
. In the literature we have found only one example where such a bound
exists. It is the case of the wavelet defined on the Fourier side by $%
\mathcal{F}\psi _{\alpha }(t)=\mathbf{1}_{\left[ 0,\infty \right] }t^{\alpha
}e^{-t}$, where the problem can be translated into a sampling problem in
Bergman spaces of analytic functions, which has been completely understood
in \cite{Seip2}.

Thanks to the connection to polyanalytic Bergman spaces, we will provide a
family of examples (which in the case $n=0$ reduces to $\psi _{\alpha }$)
where an explicit bound on the constant $b\log a$ is shown to be necessary,
providing thus the first results in this direction since the sharp
conditions associated to the Poisson wavelet. We obtain such conditions by
investigating the sampling sequences in the true polyanalytic Bergman space,
$\mathcal{A}^{n}(\mathbf{U}).$ A sampling sequence in $\mathcal{A}^{n}(%
\mathbf{U})$ is one originating inequalities of the form
\begin{equation*}
A\left\Vert F\right\Vert _{\mathcal{A}^{n}(\mathbf{U})}^{2}\leq \sum_{z\in
\Gamma (a,b)}\left\vert F(z)\right\vert ^{2}\leq B\left\Vert F\right\Vert _{%
\mathcal{A}^{n}(\mathbf{U})}^{2}.
\end{equation*}%
Thanks to the identity%
\begin{equation*}
\left\langle f\mathbf{,}T_{x}D_{s}\Phi _{n}\right\rangle _{H^{2}(\mathbf{U}%
)}=s\emph{Ber}^{n}\mathcal{F}f(z)\text{,}
\end{equation*}%
this inequality is equivalent to the wavelet frame inequality. To prove the
existence of sampling sequences, we use some results from Ascensi and Bruna
\cite{AscensiBruna}. \ Finally, we combine an argument used by Seip \cite%
{Seip} with a result from \cite{AscensiBruna} in order to prove a necessary
condition for the existence of sampling sequences in $\mathcal{A}^{n}(%
\mathbf{U})$. In terms of wavelet frames, such condiition requires that, if $%
\mathcal{W}\left( \Phi _{n},\Lambda )\right) $ is a wavelet frame for $H^{2}(%
\mathbf{U})$, then%
\begin{equation*}
b\log a<2\pi (n+1).
\end{equation*}

\subsection{Organization of the pape\textbf{r}}

The outline is as follows. We have a backgound section where we review some
known facts concerning wavelets, Laguerre functions and the fundamental
facts concerning analytic and polyanalytic Bergman spaces and the connection
between Bergman spaces and wavelets provided by the analytic Bergman
transform. Then, in the third section, we define a vector-valued version of
the continuous wavelet transform and some of its elementary properties. The
fourth section is the most important of the paper. We introduce the true
polyanalytic and the polyanalytic transforms, and describe their relation
with wavelets and vector-valued wavelets. Such a relation is fundamental in
order to prove the unitarity of the polyanalitic transforms. In section 6,
the structure of the polyanalytic Bergman spaces is investigated using these
new tools. We obtain a sequence of rational fuctions which is orthogonal in
the upper half plane, defined by its Rodrigues formula. This sequence of
rational functions is a basis of the polyanalytic Bergman space. Then we
obtain an explicit formula for the reproducing kernel of the polyanalytic
Bergman spaces, in the form of a differential operator which also ressembles
a Rodrigues formula. In the last section of the paper we present some
results concerning sampling sequences in spaces of polyanalytic functions
and their consequences in terms of wavelet frames.

\section{Background}

\subsection{The wavelet transform}

For every $x\in
\mathbb{R}
$ and $s\in
\mathbb{R}
^{+}$\ \ define the operators translation, modulation and dilation as
\begin{equation*}
T_{x}f(t)=f(t-x);\text{ \ \ \ \ \ \ \ }M_{x}f(t)=e^{-ixt}f(t),
\end{equation*}%
and
\begin{equation*}
D_{s}f(t)=s^{-\frac{1}{2}}f(s^{-1}t).
\end{equation*}%
Fix a function $g\neq 0$. Then the continuous wavelet transform of a
function $f$ with respect to a wavelet\ $g$ is defined, for every $x\in
\mathbb{R}
$, $s>0$ as
\begin{equation}
W_{g}f(x,s)=\left\langle f,T_{x}D_{s}g\right\rangle _{L^{2}(%
\mathbb{R}
)}.  \label{wavelet}
\end{equation}

The following relations are usually called \emph{the orthogonal relations
for the wavelet transform}. Assume that $g_{1},g_{2}\in L^{2}(%
\mathbb{R}
)$ satisfy

\begin{equation*}
\left\langle \mathcal{F}g_{1},\mathcal{F}g_{2}\right\rangle _{L^{2}(%
\mathbb{R}
^{+},t^{-1})}<\infty \text{.}
\end{equation*}%
Then, for all $f_{1},f_{2}\in L^{2}(%
\mathbb{R}
)$,
\begin{equation}
\left\langle W_{g_{1}}f_{1},W_{g_{2}}f_{2}\right\rangle _{L^{2}(\mathbf{U}%
,s^{-2}dxds)}=\left\langle \mathcal{F}g_{1},\mathcal{F}g_{2}\right\rangle
_{L^{2}(%
\mathbb{R}
^{+},t^{-1})}\left\langle f_{1},f_{2}\right\rangle _{L^{2}(%
\mathbb{R}
)}\text{.}  \label{ortogonalityrelations}
\end{equation}%
A function $g\in L^{2}(%
\mathbb{R}
)$ is said to be admissible if
\begin{equation*}
\left\Vert \mathcal{F}g\right\Vert _{L^{2}(%
\mathbb{R}
^{+},t^{-1})}^{2}=K\text{,}
\end{equation*}%
where $K$ is a constant. If $g$ is admissible, then for all $f\in L^{2}(%
\mathbb{R}
)$ we have
\begin{equation}
\left\Vert W_{g}f\right\Vert _{L^{2}(\mathbf{U},s^{-2}dxds)}^{2}=K\left\Vert
f\right\Vert _{L^{2}(%
\mathbb{R}
)}^{2}\text{.}  \label{isometry}
\end{equation}%
Therefore, the continuous wavelet transform provides an isometric inclusion,
being an isometry when $K=1$.%
\begin{equation*}
W_{g}:L^{2}(%
\mathbb{R}
)\rightarrow L^{2}(\mathbf{U,}s^{-2}dxds)\text{,}
\end{equation*}

If we restrict to functions $f\in H^{2}(\mathbf{U})$, the Hardy space in the
upper half plane constituted by analytic functions $f$ such that
\begin{equation*}
\text{ }\sup_{0<s<\infty }\int_{-\infty }^{\infty }\left\vert
f(x+is)\right\vert ^{2}dx<\infty \text{.}
\end{equation*}%
By the Paley-Wiener theorem,\ $H^{2}(\mathbf{U})$\ is constituted by the
functions whose Fourier transform,
\begin{equation*}
\mathcal{F}f(\omega )=\int_{%
\mathbb{R}
}f(x)e^{-ix\omega }dx,
\end{equation*}%
is supported in $%
\mathbb{R}
^{+}$ and belongs to $L^{2}(%
\mathbb{R}
^{+})$. Using the action of the Fourier transform on the dilation and
translation operators,%
\begin{equation*}
\mathcal{F}D_{s}f=D_{\frac{1}{s}}\mathcal{F}f\text{ \ \ \ \ and \ \ \ \ \ }%
\mathcal{F}T_{x}f=M_{-x}\mathcal{F}f\text{,}
\end{equation*}%
we can use Plancherel theorem to rewrite the wavelet transform
\textquotedblright on the Fourier side\textquotedblright\ as
\begin{equation*}
W_{g}f(x,s)=\left\langle \mathcal{F}f,M_{-x}D_{\frac{1}{s}}\mathcal{F}%
g\right\rangle _{L^{2}(%
\mathbb{R}
^{+})}.
\end{equation*}

\subsection{The Laguerre functions}

The Laguerre polynomials will play a central role in our discussion. One way
to define them is by the power series
\begin{equation}
L_{n}^{\alpha }(x)=\frac{(\alpha +1)_{n}}{n!}\sum_{k=0}^{n}\frac{(-n)_{k}}{%
(\alpha +1)_{k}}\frac{x^{k}}{k!}\text{.}  \label{powerseries}
\end{equation}
This is equivalent to the Rodrigues formula
\begin{equation}
L_{n}^{\alpha }(x)=\frac{e^{x}x^{-\alpha }}{n!}\frac{d^{n}}{dx^{n}}\left[
e^{-x}x^{\alpha +n}\right] \text{.}  \label{Rodr}
\end{equation}
\ The Laguerre functions are defined as
\begin{equation*}
\emph{l}_{n}^{\alpha }(x)=\mathbf{1}_{\left[ 0,\infty \right]
}(x)e^{-x/2}x^{\alpha /2}L_{n}^{\alpha }(x)\text{.}
\end{equation*}
It is well known that, for $\alpha \geq 0$, these functions constitute an
orthogonal basis for the space $L^{2}(0,\infty )$.

\section{Bergman spaces}

\subsection{Analytic and polyanalytic Bergman spaces}

With the Wirtinger diferential operator notation,%
\begin{equation*}
\frac{d}{dz}=\frac{1}{2}\left( \frac{\partial }{\partial x}-i\frac{\partial
}{\partial s}\right) \text{, \ \ }\frac{d}{d\overline{z}}=\frac{1}{2}\left(
\frac{\partial }{\partial x}+i\frac{\partial }{\partial s}\right) \text{,}
\end{equation*}%
a complex valued function $f$, is said to be analytic in a domain, if, for
every $z$ in such domain, it satisfies%
\begin{equation*}
\frac{d}{d\overline{z}}f(z)=0\text{.}
\end{equation*}%
More generally, $f$ is said to be polyanalytic of order $n$ if
\begin{equation*}
\left( \frac{d}{d\overline{z}}\right) ^{n}f(z)=0\text{.}
\end{equation*}%
Then $A(\mathbf{U})$ stands for the \emph{Bergman space} in the upper half
plane, constituted by the analytic functions in $\mathbf{U}$ such that
\begin{equation}
\int_{\mathbf{U}}\left\vert f(z)\right\vert ^{2}dxds<\infty \text{.}
\label{normbergman}
\end{equation}%
The space constituted by the polyanalytic functions of order $n$, equipped
with the same norm as the Bergman space is called the \emph{polyanalytic
Bergman space}, $\mathbf{A}^{n}(\mathbf{U})$. With this notation, $\mathbf{A}%
^{1}(\mathbf{U})=A(\mathbf{U})$. Consider also the \emph{true polyanalytic
Bergman space}, $\mathcal{A}^{n}(\mathbf{U})$, defined as
\begin{equation*}
\mathcal{A}^{n-1}(\mathbf{U})=\mathbf{A}^{n}(\mathbf{U})\ominus \mathbf{A}%
^{n-1}(\mathbf{U})\text{,}
\end{equation*}%
so that the following decomposition holds:
\begin{equation}
\mathbf{A}^{n}(\mathbf{U})=\mathcal{A}^{0}(\mathbf{U})\oplus ...\oplus
\mathcal{A}^{n-1}(\mathbf{U}).  \label{decomposition}
\end{equation}

\subsection{The Bergman transform}

We can relate the wavelet transform to Bergman spaces of analytic functions,
by choosing the window $\psi _{\alpha }$ such that
\begin{equation}
\mathcal{F}\psi _{\alpha }(t)=\mathbf{1}_{\left[ 0,\infty \right] }t^{\alpha
}e^{-t}\text{.}  \label{fg0}
\end{equation}%
Writing $z=x+si$ gives
\begin{equation}
\mathcal{F}T_{-x}D_{s}\psi _{\alpha }(t)=\mathbf{1}_{\left[ 0,\infty \right]
}s^{\alpha +\frac{1}{2}}t^{\alpha }e^{izt}  \label{Fg}
\end{equation}%
and
\begin{equation}
W_{\overline{\psi _{\alpha }}}f(-x,s)=s^{\alpha +\frac{1}{2}%
}\int_{0}^{\infty }t^{\alpha }(\mathcal{F}f)(t)e^{izt}dt\text{.}
\label{anwav}
\end{equation}%
Considering $f\in H^{2}(\mathbf{U})$, then $\mathcal{F}f\in L^{2}(%
\mathbb{R}
^{+})$. This motivates the definition of the \emph{Bergman transform }of
order\emph{\ }$\alpha $ as the analytic part of (\ref{anwav}):
\begin{equation}
Ber_{\alpha }\text{ }f(z)=s^{-\alpha }W_{\overline{\psi _{\alpha -\frac{1}{2}%
}}}(\mathcal{F}^{-1}f)(-x,s)=\int_{0}^{\infty }t^{\alpha -\frac{1}{2}%
}f(t)e^{izt}dt\text{.},  \label{Ber}
\end{equation}%
We will write
\begin{equation}
Ber\text{ }f(z)=Ber_{1}\text{ }f(z)=\int_{0}^{\infty }t^{\frac{1}{2}%
}f(t)e^{izt}dt\text{,}  \label{Beranalytic}
\end{equation}%
in order to obtain an isometric transformation
\begin{equation*}
Ber:L^{2}(%
\mathbb{R}
^{+})\rightarrow A(\mathbf{U})\text{.}
\end{equation*}

\section{A continuous vector valued wavelet transform}

Now, consider the Hilbert space $\mathcal{H}=L^{2}(%
\mathbb{R}
,\mathbf{C}^{n})$ consisting of vector-valued functions $\mathbf{f}%
=(f_{0},...,f_{n-1})$ with the inner product

\begin{equation}
\left\langle \mathbf{f,g}\right\rangle _{\mathcal{H}}=\sum_{0\leq k\leq
n-1}\left\langle f_{k},g_{k}\right\rangle _{L^{2}(%
\mathbb{R}
)}\text{.}  \label{innersuper}
\end{equation}

\begin{definition}
Let $\mathbf{g=(}g_{0},...,g_{n-1}\mathbf{)}$ be a vector of functions in $%
\mathcal{H}$ such that
\begin{equation}
\left\langle \mathcal{F}g_{i},\mathcal{F}g_{j}\right\rangle _{L^{2}(%
\mathbb{R}
^{+},t^{-1})}=\delta _{i,j}  \label{ortg}
\end{equation}%
The \emph{continuous vector valued wavelet transform} of a function $\mathbf{%
f=(}f_{1},...,f_{n-1}\mathbf{)}$ with respect to the vectorial window\ $%
\mathbf{g}$ is defined, for every $x\in
\mathbb{R}
,s\in
\mathbb{R}
^{+}$, as
\begin{equation}
\mathbf{W}_{\mathbf{g}}\mathbf{f}(x,\omega )=\left\langle \mathbf{f,}D%
\mathbf{_{s}}T_{x}\mathbf{g}\right\rangle _{\mathcal{H}}.
\label{superwavelet}
\end{equation}
\end{definition}

We can also write
\begin{equation*}
\mathbf{W}_{\mathbf{g}}\mathbf{f}(x,s)=\sum_{0\leq k\leq
n-1}W_{g_{k}}f_{k}(x,s)\text{.}
\end{equation*}
This defines a map
\begin{equation*}
\mathbf{W}_{\mathbf{g}}\mathbf{f:}\mathcal{H\rightarrow }L^{2}(\mathbf{U,}%
s^{-2}dsds)\text{.}
\end{equation*}
The orthogonality condition imposed on the vector $\mathbf{g}$ allows the
superwavelet transform to retain most of the properties of the scalar
Wavelet transform. In particular, we have vector valued versions of the
isometric property and orthogonality relations.

\begin{proposition}
Let $\mathbf{g}$ satisfy (\ref{ortg}).\ Then, for $\mathbf{f}_{1}\mathbf{,f}%
_{2},\in \mathcal{H}$,
\begin{equation}
\left\langle \mathbf{W}_{\mathbf{g}}\mathbf{f}_{1},\mathbf{W}_{\mathbf{g}}%
\mathbf{f}_{2}\right\rangle _{L^{2}(\mathbf{U,}s^{-2}dz)}=\left\langle
\mathbf{f}_{1}\mathbf{,f}_{2}\right\rangle _{\mathcal{H}}\text{.}
\label{vectorort}
\end{equation}%
In particular, $\mathbf{W}_{\mathbf{g}}\mathbf{f}$ is an isometry between
Hilbert spaces, that is
\begin{equation}
\left\Vert \mathbf{W}_{\mathbf{g}}\mathbf{f}\right\Vert _{L^{2}(\mathbf{U,}%
s^{-2}dz)}=\left\Vert \mathbf{f}\right\Vert _{\mathcal{H}}.
\label{vectorisometry}
\end{equation}
\end{proposition}

\begin{proof}
From (\ref{ortogonalityrelations}) and (\ref{ortg}),
\begin{equation}
\left\langle W_{g_{k}}f_{k},W_{g_{j}}f_{j}\right\rangle _{L^{2}(\mathbf{U,}%
s^{-2}dxds)}=\left\langle f_{k},f_{j}\right\rangle _{L^{2}(%
\mathbb{R}
)}\times \delta _{k,j}\text{.}  \label{ortsupergabor}
\end{equation}%
Then,
\begin{eqnarray*}
\left\langle \mathbf{W}_{\mathbf{g}}\mathbf{f}_{1},\mathbf{W}_{\mathbf{g}}%
\mathbf{f}_{2}\right\rangle _{L^{2}(\mathbf{U,}s^{-2}dxds)} &=&\sum_{0\leq
k,j\leq n-1}\left\langle W_{g_{k}}f_{1,k}\mathbf{,}W_{g_{j}}f_{2,k}\right%
\rangle _{L^{2}(\mathbf{U,}s^{-2}dxds)} \\
&=&\sum_{0\leq k,j\leq n-1}\left\langle f_{1,k},f_{2,j}\right\rangle _{L^{2}(%
\mathbf{U,}s^{-2}dxds)}\times \delta _{k,j} \\
&=&\sum_{0\leq k\leq n-1}\left\langle f_{1,k},f_{2,k}\right\rangle _{L^{2}(%
\mathbf{U,}s^{-2}dxds)} \\
&=&\left\langle \mathbf{f}_{1}\mathbf{,f}_{2}\right\rangle _{\mathcal{H}}%
\text{.}
\end{eqnarray*}
\end{proof}

Now let $\mathbf{W}_{\mathbf{g}}$ stand for the subspace of $L^{2}(%
\mathbb{R}
^{2})$ constituted by the image of $\mathcal{H}$ under the vector valued
wavelet transform $\mathbf{W}_{\mathbf{g}}\mathbf{f}$:
\begin{equation*}
\mathbf{W}_{\mathbf{g}}=\left\{ \mathbf{W}_{\mathbf{g}}\mathbf{f}:\mathbf{f}%
\in \mathcal{H}\right\} \text{.}
\end{equation*}%
Since
\begin{equation*}
\mathbf{W}_{\mathbf{g}}\mathbf{f}=\sum_{0\leq k\leq n-1}W_{g_{k}}f_{k}
\end{equation*}%
and
\begin{equation*}
\left\langle W_{g_{k}}f_{k},W_{g_{j}}f_{j}\right\rangle _{L^{2}(\mathbf{U,}%
s^{-2}dxds)}=\delta _{k,j}\text{,}
\end{equation*}%
we know that every $F\in \mathbf{W}_{g}$ can be written in a unique way in
the form
\begin{equation}
F=F_{0}+...+F_{n-1}\text{.}  \label{superF}
\end{equation}%
As a result,
\begin{equation}
\mathbf{W}_{\mathbf{g}}=\mathcal{W}_{g_{0}}\oplus ...\oplus \mathcal{W}%
_{g_{n-1}}\text{,}  \label{Superdecomposition}
\end{equation}%
where%
\begin{equation*}
\mathcal{W}_{g_{j}}=\{W_{g_{j}}f:f\in L^{2}(%
\mathbb{R}
)\}.
\end{equation*}

\begin{proposition}
The space $\mathbf{W}_{\mathbf{g}}$ is a Hilbert space with reproducing
kernel given by
\begin{equation}
\mathbf{k}(z,w)=\left\langle T_{\eta }D_{u}\mathbf{g},T_{x}D_{s}\mathbf{g}%
\right\rangle _{\mathcal{H}}=\sum_{j=0}^{n-1}k_{j}\left( z,w\right) \text{,}
\label{superrep}
\end{equation}%
where $k_{j}\left( z,w\right) $ is the reproducing kernel of $\mathcal{W}%
_{g_{j}}$.
\end{proposition}

\begin{proof}
Let $\mathbf{F}\in \mathbf{W}_{\mathbf{g}}$. There exists $\mathbf{f}\in
\mathcal{H}$ such that $\mathbf{F=W}_{\mathbf{g}}\mathbf{f}$. By definition,
$\mathbf{k}(z,.)=\mathbf{W}_{g}\mathbf{(}T_{x}D\mathbf{_{s}\mathbf{g})}$.
Thus, using (\ref{vectorort}),
\begin{eqnarray*}
\left\langle \mathbf{F},\mathbf{k(}z,.\mathbf{)}\right\rangle _{L^{2}(%
\mathbf{U,}s^{-2}dz)} &=&\left\langle \mathbf{W}_{\mathbf{g}}\mathbf{f},%
\mathbf{W}_{g}\mathbf{(}T_{x}D_{s}\mathbf{\mathbf{g})}\right\rangle _{L^{2}(%
\mathbf{U,}s^{-2}dz)} \\
&=&\left\langle \mathbf{f,}T_{x}D_{s}\mathbf{\mathbf{g}}\right\rangle _{%
\mathcal{H}} \\
&=&\mathbf{F(}z\mathbf{)}.
\end{eqnarray*}%
The second inequality follows from the well known fact that the reproducing
kernel of the space $\mathcal{W}_{g_{j}}$ is given by $\left\langle T_{\eta
}D_{u}g_{j},T_{x}D_{s}g_{j}\right\rangle _{L^{2}(%
\mathbb{R}
)}.$
\end{proof}

Through the paper, we will restrict ourselfs to vectors $\mathbf{f}$ such
that the Fourier transform of each of its components is supported in $%
\mathbb{R}
^{+}$ and belongs to $L^{2}(%
\mathbb{R}
^{+})$. In such a case, $\mathcal{H}=H^{2}(\mathbf{U},\mathbf{C}^{n})$ and,
on the Fourier side, the notation $\mathcal{H}^{+}=L^{2}(%
\mathbb{R}
^{+},\mathbf{C}^{n})$, will be used.

\section{The polyanalytic Bergman transform}

In this section we will study a special case of the continuous vector valued
wavelet transform, when the vector is defined in terms of Laguerre
functions. First we treat the scalar case, which originates a unitary map
onto the \emph{true} polyanalytic Bergman space.We show that the required
unitary mappings can be related to special wavelet transforms and, via a
connection to the previous section, we define the vector valued polyanalytic
transform onto the polyanalytic space (which we call the polyanalytic
Bergman transform).

\subsection{The true polyanalytic Bergman transform}

First we will study the transform that allows, in the multiplexing context
explained in the introduction, to send each signal $f_{k}\in L^{2}(%
\mathbb{R}
^{+})$ to a space $\mathcal{A}^{k}(\mathbf{U}).$

\begin{definition}
The \emph{true polyanalytic Bergman transform} of order $n$ is the transform
mapping every $f\in L^{2}(%
\mathbb{R}
^{+})$ to
\begin{equation}
\emph{Ber}^{n}f(z)=\frac{1}{(2i)^{n}n!}\left( \frac{d}{dz}\right) ^{n}\left[
s^{n}F(z)\right] \text{,}  \label{polyber}
\end{equation}%
where $F=Ber$ $f$ and $f\in L^{2}(%
\mathbb{R}
^{+})$.
\end{definition}

The purpose of this section is to prove that the transform $\emph{Ber}^{n}$
is unitary%
\begin{equation*}
\emph{Ber}^{n}:L^{2}(%
\mathbb{R}
^{+})\rightarrow \mathcal{A}^{n}(\mathbf{U})\text{.}
\end{equation*}
We will need some identities which have independent interest. First observe
that, since
\begin{equation*}
L_{n}^{0}(t)=\sum_{k=0}^{n}(-1)^{k}\binom{n}{k}\frac{t^{k}}{k!}\text{,}
\end{equation*}%
we have
\begin{equation}
\emph{l}_{n}^{0}(t)=\sum_{k=0}^{n}(-1)^{k}\binom{n}{k}\frac{1}{k!}\mathbf{1}%
_{\left[ 0,\infty \right] }(x)t^{k}e^{-t/2}=t^{-\frac{1}{2}%
}\sum_{k=0}^{n}(-1)^{k}\binom{n}{k}\mathcal{F}\psi _{k+\frac{1}{2}}(\frac{t}{%
2})\text{.}  \label{comb}
\end{equation}%
Due to this observation, it is reasonable to expect that the functions $\Phi
_{n}$, defined by
\begin{equation}
\mathcal{F}\Phi _{n}(t)=t^{\frac{1}{2}}\emph{l}_{n}^{0}(2t)\text{,}
\label{defS}
\end{equation}%
will play a distinguished role in our analysis. Indeed, an essential step in
the proof of the unitary property is to write (\ref{polyber}) in terms of a
wavelet transform with analysing wavelet $\Phi _{n}$.

\begin{proposition}
The true polyanalytic Bergman transform of order $n$ can be written as:
\end{proposition}

\begin{enumerate}
\item A polyanalytic function of order $n+1$:
\begin{equation}
\emph{Ber}^{n}f(z)=\frac{1}{(2i)^{n}n!}\sum_{k=0}^{n}(2i)^{k}\binom{n}{k}%
\frac{1}{k!}s^{k}F^{(k)}(z).  \label{poly}
\end{equation}

\item In terms of analytic Bergman transforms of different orders:
\begin{equation}
\emph{Ber}^{n}f(z)=\sum_{k=0}^{n}(-2)^{k}\binom{n}{k}s^{k}Ber_{k+1}\text{ }%
f(z)\text{.}  \label{orders}
\end{equation}

\item In terms of a wavelet transform:
\begin{equation}
\emph{Ber}^{n}f(z)=s^{-1}W_{\Phi _{n}}(\mathcal{F}^{-1}f)(x,s)\text{.}
\label{waveberg}
\end{equation}
\end{enumerate}

\begin{proof}
The first identity follows from a standard application of Leibnitz formula.
Then, since
\begin{equation*}
\frac{d}{d\overline{z}}s^{k}=\frac{1}{2}ks^{k-1}\text{,}
\end{equation*}%
we have
\begin{equation*}
\left( \frac{d}{d\overline{z}}\right) ^{n+1}\emph{Ber}^{n}f(z)=0\text{,}
\end{equation*}%
and $\emph{Ber}^{n}f$ is polyanalytic of order $n+1$. To prove (\ref{orders}%
), observe that diferentiating (\ref{Beranalytic}) under the integral sign
gives
\begin{equation}
F^{(k)}(z)=\left( \frac{d}{dz}\right) ^{k}Ber\text{ }f(z)=i^{k}Ber_{k+1}%
\text{ }f(z)\text{.}  \label{derivaBergman}
\end{equation}%
Applying this to (\ref{poly}) gives (\ref{orders}).

Now (\ref{waveberg}). Combining (\ref{comb}) with (\ref{defS}) and inverting
the Fourier transform gives
\begin{equation*}
\Phi _{n}(t)=\sum_{k=0}^{n}(-2)^{k}\binom{n}{k}\psi _{k+\frac{1}{2}}(t)\text{%
.}
\end{equation*}%
Then,
\begin{eqnarray*}
W_{\Phi _{n}}f(x,s) &=&\sum_{k=0}^{n}(-2)^{k}\binom{n}{k}W_{\overline{\psi
_{k+\frac{1}{2}}}}f(-x,s) \\
&=&\sum_{k=0}^{n}(-2)^{k}\binom{n}{k}s^{k+1}Ber_{k+1}\text{ }\mathcal{F}f(z)
\\
&=&s\emph{Ber}^{n}(\mathcal{F}f)(z)\text{.}
\end{eqnarray*}%
and (\ref{waveberg}) follows.
\end{proof}

Now we can prove the main result. The idea consists in writing the wavelet
transform (\ref{waveberg}) as a composition of several unitary operators and
is suggested by the techniques used in \cite{VasiBergman}\ and \cite{Hutnikm}%
. We need to introduce two auxiliary operators. For convenience write $L^{2}(%
\mathbf{U})=L^{2}(%
\mathbb{R}
)\mathcal{\otimes }L^{2}(%
\mathbb{R}
^{+})$ and define the unitary operators $U_{1,2}:L^{2}(%
\mathbb{R}
)\mathcal{\otimes }L^{2}(%
\mathbb{R}
^{+})\rightarrow L^{2}(%
\mathbb{R}
)\mathcal{\otimes }L^{2}(%
\mathbb{R}
^{+})$:
\begin{eqnarray*}
U_{1}(F)(x,s) &=&(\mathcal{F}^{-1}\mathcal{\otimes }I)(F\mathcal{)}(x,s) \\
U_{2}(F)(x,s) &=&\frac{1}{\sqrt{2\left\vert x\right\vert }}F(x,\frac{s}{%
2\left\vert x\right\vert })\text{.}
\end{eqnarray*}

We will need the following result of Vasilevski \textbf{\cite{VasiBergman}}.

\textbf{Theorem A \cite[Corollary 4.2]{VasiBergman} }\emph{Let }$L_{n}$\emph{%
\ stand for the space generated by }$1_{\left[ 0,\infty \right] }l_{n}^{0}$%
\emph{\ \ The operator }$U=U_{2}U_{1}$,%
\begin{equation*}
U:\mathcal{A}^{n}(\mathbf{U})\rightarrow L^{2}(%
\mathbb{R}
^{+})\mathcal{\otimes }\emph{L}_{n}
\end{equation*}%
such that, given $f\in \mathcal{A}^{n}(\mathbf{U})$,%
\begin{equation*}
(Uf)(x,s)=\mathbf{1}_{\left[ 0,\infty \right] }(x)f(x)\emph{l}_{n}^{0}(s)%
\text{,}
\end{equation*}%
\emph{is unitary.}

We now combine Theorem A with Proposition 3 to prove the main result.

\begin{theorem}
The transform $\emph{Ber}^{n}:L^{2}(%
\mathbb{R}
^{+})\rightarrow \mathcal{A}^{n}(\mathbf{U})$ is unitary.
\end{theorem}

\begin{proof}
Consider the unitary operator
\begin{equation*}
R_{k}:L^{2}(%
\mathbb{R}
^{+})\rightarrow L^{2}(%
\mathbb{R}
^{+})\mathcal{\otimes }\emph{L}_{k}
\end{equation*}%
defined by
\begin{equation*}
(R_{k}f)(x,s)=\mathbf{1}_{\left[ 0,\infty \right] }(x)f(x)\emph{l}_{k}^{0}(s)%
\text{.}
\end{equation*}%
Then the composition \ $U^{-1}R_{k}$ is also unitary
\begin{equation*}
U^{-1}R_{k}:L^{2}(%
\mathbb{R}
^{+})\rightarrow \mathcal{A}^{k}(\mathbf{U})
\end{equation*}%
We will now show that this transform is exactly $\emph{Ber}^{n}$. From the
definition of $U_{2}$ it is easy to see that
\begin{equation*}
U_{2}^{-1}(F)(x,s)=\sqrt{2\left\vert x\right\vert }F(x,2\left\vert
x\right\vert s)
\end{equation*}%
and
\begin{equation*}
(U_{2}^{-1}R_{k}f)(x,s)=\mathbf{1}_{\left[ 0,\infty \right] }(x)\sqrt{2x}f(x)%
\emph{l}_{k}^{0}(2xs).
\end{equation*}%
Applying $U_{1}^{-1}$ gives
\begin{eqnarray*}
(U^{-1}R_{k}f)(x,s) &=&s^{-1}\int_{0}^{\infty }f(t)s^{\frac{1}{2}}(2ts)^{%
\frac{1}{2}}\emph{l}_{k}^{0}(2ts)e^{ixt}dt \\
&=&s^{-1}\int_{%
\mathbb{R}
}(\mathcal{F}^{-1}f)(t)s^{-\frac{1}{2}}\overline{\Phi _{n}(s^{-1}(t-x))}dt \\
&=&s^{-1}W_{\Phi _{n}}(\mathcal{F}^{-1}f)(x,s) \\
&=&\emph{Ber}^{n}f(z)\text{,}
\end{eqnarray*}%
by using the identity (\ref{waveberg}).
\end{proof}

\subsection{The polyanalytic Bergman transform}

Now, consider the Hilbert space $\mathcal{H}^{+}=L^{2}(%
\mathbb{R}
^{+},\mathbf{C}^{n})$ consisting of vector-valued functions $\mathbf{f}%
=(f_{0},...,f_{n-1})$ with the inner product%
\begin{equation*}
\left\langle \mathbf{f,g}\right\rangle _{\mathcal{H}^{+}}=\sum_{0\leq k\leq
n-1}\left\langle f_{k},g_{k}\right\rangle _{L^{2}(%
\mathbb{R}
^{+})}\text{.}
\end{equation*}

\begin{definition}
The \emph{polyanalytic Bergman transform of order n} is defined, for $%
\mathbf{f}\in \mathcal{H}^{+}$ as
\begin{equation*}
\mathbf{Be}r^{n}\mathbf{f}=\sum_{0\leq k\leq n-1}Ber^{k}f_{k}\text{.}
\end{equation*}
\end{definition}

If we take, then we have the following relation with the polyanalytic
Bergman transform:

\begin{theorem}
The \emph{polyanalytic Bergman transform of order n} is a unitary operator
\begin{equation*}
\mathbf{Be}r^{n}:\mathcal{H}^{+}\rightarrow \mathbf{A}^{n}(\mathbf{U})
\end{equation*}
\end{theorem}

\begin{proof}
To see that it is onto, let $\mathbf{F}\in \mathbf{A}^{n}(\mathbf{U})$.
Then, using (\ref{decomposition}), write%
\begin{equation*}
\mathbf{F=}F_{0}+...+F_{n-1},
\end{equation*}%
with $F_{k}\in \emph{A}^{k}(\mathbf{U}),$ $k=0,...,n-1.$ Since $Ber^{k}$ is
onto, for every $k=0,...,n-1$ there exists $f_{k}\in L^{2}(%
\mathbb{R}
^{+})$ such that $F_{k}=Ber^{k}f_{k}.$ To prove the isometry, we first
relate the polyanalytic Bergman transform to the vector valued wavelet
transform with the vectorial window $\mathbf{\Phi }_{n}=(\Phi _{0},...,\Phi
_{n-1})$, using the identity $\emph{Ber}^{n}f(z)=s^{-1}W_{\Phi _{n}}(%
\mathcal{F}^{-1}f)(x,s)$:
\begin{equation*}
\mathbf{W}_{\mathbf{\Phi }_{n}}(\mathcal{F}^{-1}\mathbf{f)}(x,s)=\sum_{0\leq
k\leq n-1}W_{\mathbf{\Phi }_{k}}(\mathcal{F}^{-1}f_{k})(x,s)=\sum_{0\leq
k\leq n-1}sBer^{k}f_{k}=s\mathbf{Ber}^{n}\mathbf{f}\text{.}
\end{equation*}%
Now, combining this with (\ref{vectorisometry}),%
\begin{equation*}
\left\Vert \mathbf{Ber}^{n}\mathbf{f}\right\Vert _{\mathbf{A}^{n}(\mathbf{U}%
)}=\left\Vert \mathbf{W}_{\mathbf{\Phi }_{n}}(\mathcal{F}^{-1}\mathbf{f)}%
\right\Vert _{\mathbf{L}^{2}(\mathbf{U,}s^{-2}dxds)}=\left\Vert \mathcal{F}%
^{-1}\mathbf{f}\right\Vert _{\mathcal{H}}=\left\Vert \mathbf{f}\right\Vert _{%
\mathcal{H}^{+}}\text{.}
\end{equation*}
\end{proof}

\section{The structure of polyanalytic Bergman spaces}

The purpose of this section is to apply the connection to wavelet transforms
to study polyanalytic Bergman spaces. We will obtain an orthogonal basis for
$\mathbf{A}^{n}(\mathbf{U})$ and compute an explicit formula for the
reproducing kernel. The reproducing kernel, $K^{n}(z,w)$ , of the true
polyanalytic Bergman space $\mathcal{A}^{n}(\mathbf{U})$ is also very
important, since once we have a function $F\in \mathbf{A}^{n}(\mathbf{U})$,
we can recover its true polyanalytic component $F_{k}\in \mathcal{A}^{k}(%
\mathbf{U})$ by the orthogonal projecion over the space $\mathcal{A}^{k}(%
\mathbf{U})$, which is given by the formula%
\begin{equation*}
F_{k}(z)=\left\langle F(w),K^{k}(z,w)\right\rangle _{\mathbf{A}^{n}(\mathbf{U%
})}.
\end{equation*}%
Our formulas will be given in a form of differential operators which are
reminiscent of the \emph{Rodrigues formula, }a well known structure formula
in the theory of classic orthogonal polynomials.

\subsection{An orthogonal basis}

Consider the functions $\Psi _{n}^{\beta }$ , for every $n\geq 0$ and $\beta
>1$:
\begin{equation*}
\Psi _{n}^{\beta }(z)=(2i)^{\beta +1}\frac{\Gamma (\beta +n)}{n!}\left(
\frac{z-i}{z+i}\right) ^{n}\left( \frac{1}{z+i}\right) ^{\beta }\text{.}
\end{equation*}%
It is well known that these functions constitute a basis of $A_{\beta -2}(%
\mathbf{U})$. A calculation using the special function formula
\begin{equation*}
\int_{0}^{\infty }x^{\alpha }L_{n}^{\alpha }(x)e^{-xs}dx=\frac{\Gamma
(\alpha +n+1)}{n!}s^{-\alpha -n-1}(s-1)^{n}
\end{equation*}%
gives
\begin{equation*}
Ber_{\alpha }\text{ }l_{n}^{2\alpha -1}=\Psi _{n}^{2\alpha }\text{.}
\end{equation*}%
Now write
\begin{equation*}
\Psi _{n}(z)=\Psi _{n}^{2}(z)=(2i)^{3}\frac{\Gamma (2+n)}{n!}\left( \frac{z-i%
}{z+i}\right) ^{n}\left( \frac{1}{z+i}\right) ^{2}
\end{equation*}%
to denote a basis of $A(\mathbf{U})$, so that
\begin{equation*}
Ber\text{ }l_{n}^{1}=\Psi _{n}\text{.}
\end{equation*}

\begin{definition}
Define a set of functions by
\begin{equation*}
e_{n,m}(z)=\frac{1}{(2i)^{n}n!}\left( \frac{d}{dz}\right) ^{n}\left[
s^{n}\Psi _{m}(z)\right] \text{.}
\end{equation*}
\end{definition}

\begin{proposition}
The set $\{e_{k,m}\}_{k\geq 0,0\leq m<n}$ is an orthonormal basis of $%
\mathbf{A}^{n}(\mathbf{U})$.
\end{proposition}

\begin{proof}
Since%
\begin{equation*}
e_{n,k}(z)=\emph{Ber}^{n}l_{k}^{1}=s^{-1}W_{\Phi _{n}}(\mathcal{F}^{-1}\emph{%
l}_{k}^{1})(x,s),
\end{equation*}%
the orthogonality follows from (\ref{ortogonalityrelations}):
\begin{eqnarray*}
\left\langle e_{n,k},e_{l,j}\right\rangle _{L^{2}(\mathbf{U},dz)}
&=&\left\langle W_{\Phi _{n}}(\mathcal{F}^{-1}\emph{l}_{k}^{1}),W_{\Phi
_{l}}(\mathcal{F}^{-1}\emph{l}_{j}^{1})\right\rangle _{L^{2}(\mathbf{U}%
,s^{-2}dz)} \\
&=&\left\langle \emph{l}_{n}^{0},\emph{l}_{l}^{0}\right\rangle _{L^{2}(%
\mathbb{R}
^{+})}\left\langle \mathcal{F}^{-1}\emph{l}_{k}^{1},\mathcal{F}^{-1}\emph{l}%
_{j}^{1}\right\rangle _{H^{2}(\mathbf{U})} \\
&=&\delta _{n,l}\delta _{k,j}\text{.}
\end{eqnarray*}%
The unitarity of $\emph{Ber}^{n}$ shows that, for every $m$, $%
\{e_{k,m}\}_{k\geq 0}$ spans $\mathcal{A}^{m}(\mathbf{U})$, since $\{\emph{l}%
_{n}^{1}\}_{n\geq 0}$ spans $L^{2}(%
\mathbb{R}
^{+})$. From the decomposition (\ref{decomposition}), every element in $%
\mathbf{A}^{n}(\mathbf{U})$\ can be written as a linear combination of
elements of $\{\mathcal{A}^{m}(\mathbf{U})\}_{m<n}$.\ Therefore $%
\{e_{k,m}\}_{k\geq 0,0\leq m<n}$ spans $\mathbf{A}^{n}(\mathbf{U})$.
\end{proof}

\begin{corollary}
The set $\{e_{k,m}\}_{,0\leq m<n}$ is an orthonormal basis of $\mathcal{A}%
^{k}(\mathbf{U})$.
\end{corollary}

\begin{proof}
This follows immediately from the decomposition (\ref{decomposition}).
\end{proof}

\subsection{The reproducing kernel}

In our computations of the reproducing kernels we will need the following
relation, which says essentially that the Bergman transform intertwines the
representation of the affine group in $L^{2}(%
\mathbb{R}
^{+})$ with its representation in the Bergman space:
\begin{equation}
\emph{Ber}(M_{-\mu }D_{\frac{1}{\eta }}f)=s^{-1}\emph{Ber}(f)(\frac{z+u}{%
\eta })\text{;}  \label{intertwining}
\end{equation}%
The identity (\ref{intertwining}) follows from the change of variables:%
\begin{equation*}
\int_{0}^{\infty }t^{\frac{1}{2}}e^{izt}(M_{-\mu }D_{\frac{1}{\eta }%
}f)(t)dt=\int_{0}^{\infty }t^{\frac{1}{2}}e^{i(\frac{z+u}{\eta })t}f(t)dt
\end{equation*}%
It is interesting to observe that, defining a transform by $T=\emph{Ber}(%
\mathcal{F}f)$ and $Mf(t)=t^{\frac{1}{2}}f(t)$, we have the following
commutative diagram:%
\begin{equation*}
\begin{array}{ccc}
H^{2}(\mathbf{U}) & _{\longrightarrow }^{T} & A^{2}(\mathbf{U}) \\
\mathcal{F}\downarrow &  & \downarrow \mathcal{F} \\
L^{2}(%
\mathbb{R}
^{+}) & _{\longrightarrow }^{M} & L^{2}(%
\mathbb{R}
^{+},t^{-1}dt)%
\end{array}%
\end{equation*}%
The Fourier isometry on the right column is the Paley-Wiener theorem for the
Bergman space \cite{DGM}.

Now we compute the reproducing kernel of the wavelet space $\mathcal{W}%
_{\Phi _{n}}$.

\begin{theorem}
The reproducing kernel of $\mathcal{W}_{\Phi _{n}}$is given by
\begin{equation*}
k^{n}(z,w)=\frac{1}{n!(2i)^{n}}\frac{s}{\eta }\left( \frac{d}{dz}\right) ^{n}%
\left[ s^{n}\Omega _{n}\left( \frac{z+u}{\eta }\right) \right] ,
\end{equation*}%
where%
\begin{equation*}
\Omega _{n}(z)=4(n+z-i)\left( \frac{1}{z+i}\right) ^{3}\left( \frac{z-i}{z+i}%
\right) ^{n-1}.
\end{equation*}
\end{theorem}

\begin{proof}
The reproducing kernel of $\mathcal{W}_{\Phi _{n}}$ is
\begin{eqnarray*}
k^{n}(z,w) &=&\left\langle T_{-\mu }D_{\eta }\Phi _{n},T_{-x}D_{s}\Phi
_{n}\right\rangle _{H^{2}(\mathbf{U})} \\
&=&s\emph{Ber}^{n}(M_{-\mu }D_{\frac{1}{\eta }}\mathcal{F}\Phi _{n})(z) \\
&=&\frac{s}{n!(2i)^{n}}\sum_{k=0}^{n}(2i)^{k}\binom{n}{k}\frac{s^{k}}{k!}%
\left( \frac{d}{dz}\right) ^{k}\left[ \emph{Ber}(M_{-\mu }D_{\frac{1}{\eta }}%
\mathcal{F}\Phi _{n})(z)\right]  \\
&=&\frac{s}{n!(2i)^{n}}\left( \frac{d}{dz}\right) ^{n}\left[ s^{n}\emph{Ber}%
(M_{-\mu }D_{\frac{1}{\eta }}\mathcal{F}\Phi _{n})(z)\right]
\end{eqnarray*}%
Now, (\ref{intertwining}) gives%
\begin{equation*}
k^{n}(z,w)=\frac{1}{(2i)^{n}n!}\frac{s}{\eta }\left( \frac{d}{dz}\right) ^{n}%
\left[ s^{n}\emph{Ber}(\mathcal{F}\Phi _{n})\left( \frac{z+u}{\eta }\right) %
\right] .
\end{equation*}%
We just need to compute%
\begin{eqnarray*}
\emph{Ber}(\mathcal{F}\Phi _{n}) &=&\int_{0}^{\infty
}tl_{n}^{0}(2t)e^{itz}dt=\frac{1}{i}\frac{d}{dz}\int_{0}^{\infty
}l_{n}^{0}(2t)e^{itz}dt \\
&=&\frac{4}{i}\frac{d}{dz}\left[ \left( \frac{z-i}{z+i}\right) ^{n}\frac{1}{%
z+i}\right] =\Omega _{n}(z)\text{,}
\end{eqnarray*}%
and the formula is proved.
\end{proof}

The next Lemma can be used to transfer properties from the spaces $\mathcal{W%
}_{\Phi _{n}}$ to the spaces $\mathcal{A}^{n}(\mathbf{U}).$

\begin{lemma}
The operator%
\begin{equation*}
E:\mathcal{W}_{\Phi _{n}}\rightarrow \mathcal{A}^{n}(\mathbf{U})
\end{equation*}%
\begin{equation*}
f\rightarrow s^{-1}f(x,s)\text{,}
\end{equation*}%
is unitary.
\end{lemma}

\begin{proof}
Clearly, $E$ is isometric. Since $l_{m}^{1}$ is a basis of $L^{2}(%
\mathbb{R}
^{+})$, then $W_{\Phi _{n}}(\mathcal{F}^{-1}\emph{l}_{m}^{1})$ is a basis of
$\mathcal{W}_{\Phi _{n}}$. Then
\begin{equation*}
E\left( W_{\Phi _{n}}(\mathcal{F}^{-1}\emph{l}_{m}^{1})\right)
=s^{-1}W_{\Phi _{n}}(\mathcal{F}^{-1}\emph{l}_{m}^{1})(x,s)=\emph{Ber}%
^{n}l_{m}^{1}(z)=e_{n,m}(z)\text{.}
\end{equation*}%
Thus, $E\left( \mathcal{W}_{\Phi _{n}}\right) $ is dense in $\mathcal{A}^{n}(%
\mathbf{U})$.
\end{proof}

\begin{theorem}
The reproducing kernels of the spaces$\ \mathcal{A}^{n}(\mathbf{U})$, $%
K^{n}(z,w)$, are given by
\begin{equation*}
K^{n}(z,w)=\frac{1}{n!(2i)^{n}}\left( \frac{d}{dz}\right) ^{n}\left[
s^{n}\Omega _{n}\left( \frac{z+u}{\eta }\right) \right] ,
\end{equation*}%
The reproducing kernels of the spaces, $\mathbf{A}^{n}(\mathbf{U})$, $%
\mathbf{K}^{n}(z,w)$, are given by
\begin{equation*}
\mathbf{K}^{n}(z,w)=\sum_{k=0}^{n-1}\frac{1}{n!(2i)^{n}}\left( \frac{d}{dz}%
\right) ^{n}\left[ s^{n}\Omega _{n}\left( \frac{z+u}{\eta }\right) \right] .
\end{equation*}
\end{theorem}

\begin{proof}
Let $f\in \mathcal{A}^{n}(\mathbf{U})$. Then, by the above Lemma, $sf(z)\in
\mathcal{W}_{\Phi _{n}}$. Therefore,
\begin{equation*}
sf(z)=\left\langle k^{n}(z,w),\eta f(w)\right\rangle _{\mathcal{W}_{\Phi
_{n}}}\text{,}
\end{equation*}%
and%
\begin{equation*}
f(z)=\left\langle \frac{\eta }{s}k^{n}(z,w),f(w)\right\rangle _{\mathbf{A}%
^{n}(\mathbf{U})}.
\end{equation*}%
We conclude that $K^{n}(z,w)=\frac{\eta }{s}k^{n}(z,w)$. The second
assertion follows imediately from (\ref{superrep}).
\end{proof}

\section{Sampling sequences and wavelet frames}

This section is devoted to sampling and frames for Wavelet frames and
super-frames. There are several approaches to general sampling and stability
problems (see, for instance \cite{FG} and \cite{WsunHans}), but we will
follow mainly the one in \cite{AscensiBruna}, where, fixed an analyzing
wavelet, the space of all continuous transforms (the "model space") is
considered, in order to translate the frame problem in a sampling problem
for such a model space.

Now, we will denote by $\Gamma (a,b)$ the set $\{z_{mk}=a^{m}(bk+i)\}.$ We
say that $\Gamma $ is \emph{a} \emph{sampling sequence} for $\mathcal{A}^{n}(%
\mathbf{U})$\ if there exist $A,B>0$ such that, for every $F\in \mathcal{A}%
^{n}(\mathbf{U})$,%
\begin{equation}
A\left\Vert F\right\Vert _{\mathcal{A}^{n}(\mathbf{U})}^{2}\leq \sum_{z\in
\Gamma (a,b)}s^{2}\left\vert F(z)\right\vert ^{2}\leq B\left\Vert
F\right\Vert _{\mathcal{A}^{n}(\mathbf{U})}^{2}.  \label{samp}
\end{equation}

We say that $\mathcal{W}\left( \psi ,\Gamma (a,b))\right) $ is a wavelet
frame for $H^{2}(\mathbf{U})$ if
\begin{equation}
A\left\Vert f\right\Vert _{H^{2}(\mathbf{U})}^{2}\leq \sum_{j,k}\left\vert
\left\langle f,T_{a^{j}bk}D_{a^{j}}\psi \right\rangle \right\vert ^{2}\leq
B\left\Vert f\right\Vert _{H^{2}(\mathbf{U})}^{2}\text{.}  \label{waveframe}
\end{equation}%
Since
\begin{equation*}
\left\langle f\mathbf{,}T_{x}D_{s}\Phi _{n}\right\rangle _{H^{2}(\mathbf{U}%
)}=W_{\Phi _{n}}f(x,s)=s\emph{Ber}^{n}\mathcal{F}f(z)\text{,}
\end{equation*}%
it is plain that $\Gamma $ is \emph{a} \emph{sampling sequence} for $%
\mathcal{A}^{n}(\mathbf{U})$\ if and only if $\mathcal{W}\left( \psi ,\Gamma
(a,b))\right) $ is a wavelet frame for $H^{2}(\mathbf{U})$. Our next result
is an upper bound on the size of the parameters $(a,b)$ (or a lower bound on
density) necessary to generate sampling sequences in the true polyanalytic
Bergman space, or wavelet frames with Laguerre functions.

The vector valued system $\mathcal{W}(\mathbf{g},\Lambda
)=\{T_{a^{j}bk}D_{a^{j}}\mathbf{g}\}$ is a $\emph{wavelet}$ \emph{superframe}
for $\mathcal{H}$ if there exist constants $A$ and $B$ such that, for every $%
\mathbf{f}\in \mathcal{H}$,
\begin{equation}
A\left\Vert \mathbf{f}\right\Vert _{\mathcal{H}}^{2}\leq
\sum_{j,k}\left\vert \left\langle \mathbf{f},T_{a^{j}bk}D_{a^{j}}\mathbf{g}%
\right\rangle _{\mathcal{H}}\right\vert ^{2}\leq B\left\Vert \mathbf{f}%
\right\Vert _{\mathcal{H}}^{2}.  \label{superframe}
\end{equation}%
Superframes were introduced in a more abstract form in \cite{DL} and in the
context of "multiplexing" in \cite{Balan}. Take the analyzing vector to be $%
\mathbf{\Phi }_{n}=(\Phi _{0},...,\Phi _{n-1})$\ Using the identity%
\begin{equation*}
\left\langle \mathbf{f,}D\mathbf{_{s}}T_{x}\mathbf{\Phi }_{n}\right\rangle _{%
\mathcal{H}}=s\mathbf{Ber}^{n}\mathbf{f}\text{,}
\end{equation*}%
we see that $\mathcal{W}(\mathbf{\Phi }_{n},\Lambda )$ is a wavelet
superframe for $\mathcal{H}$\ if and only if $\Gamma $ is a sampling
sequence for $\mathbf{A}^{n}(\mathbf{U})$.

\subsection{Existence of sampling sequences and frames}

In this subsection we will prove the existence of wavelet frames with the
functions $\Phi _{n}$, provided the hyperbolic lattice is sufficiently
dense. We need some notations from \textbf{\cite{AscensiBruna}. }First
recall the hyperbolic distance in the half-plane%
\begin{equation*}
d(z_{1},z_{2})=\frac{1}{2}\log \frac{1+\rho (z_{1},z_{2})}{1-\rho
(z_{1},z_{2})}\text{,}
\end{equation*}%
where $\rho (z_{1},z_{2})$ is the pseudohyperbolic distance:%
\begin{equation*}
\rho (z_{1},z_{2})=\left\vert \frac{z_{1}-z_{2}}{z_{1}-\overline{z_{2}}}%
\right\vert \text{.}
\end{equation*}%
Given a continuous function $h$ in $%
\mathbb{R}
\times
\mathbb{R}
^{+}$, we define its local maximal function as%
\begin{equation*}
Mh(z)=\sup_{w\in B(z,1)}\left\vert h(w)\right\vert \text{,}
\end{equation*}%
where $B(z,1)$ is the hyperbolic ball of center $z$ and radius $1$ in $%
\mathbb{R}
\times
\mathbb{R}
^{+}$. Now, define
\begin{equation*}
k_{g}(z)=\left\langle \psi \mathbf{,}T_{x}D_{s}g\right\rangle _{H^{2}(%
\mathbf{U})}\text{, }z=x+is\text{,}
\end{equation*}%
and%
\begin{equation*}
MB=\{g:Mk_{g}\in L^{1}(\mathbf{U})\}
\end{equation*}%
With these notations, Theorem 4.9 in \cite{AscensiBruna} reads:

\textbf{Theorem B. }\emph{Let }$W_{g}$\emph{\ be the wavelet space
associated to the analyzing wavelet }$g$\emph{. If }$g\in MB$\emph{, then
there is a }$\delta $\emph{\ such that every uniformly discrete set }$\Gamma
$\emph{\ satisfying }$B(z,\delta )\cap \Gamma \neq \varnothing $\emph{\ for
every }$z$\emph{, is a sampling set for }$W_{g}$\emph{.}

\begin{lemma}
$\Phi _{n}(z)\in MB$.
\end{lemma}

\begin{proof}
Using (\ref{comb}) twice and evaluating the resulting wavelet transforms on
the Fourier side, gives:%
\begin{eqnarray*}
\left\langle \Phi _{n}\mathbf{,}T_{x}D_{s}\Phi _{n}\right\rangle _{H^{2}(%
\mathbf{U})} &=&\sum_{k=0}^{n}\sum_{j=0}^{n}(-2)^{k+j}\binom{n}{k}\binom{n}{j%
}\left\langle \psi _{j+\frac{1}{2}},T_{x}D_{s}\psi _{k+\frac{1}{2}%
}\right\rangle _{H^{2}(\mathbf{U})} \\
&=&\sum_{k=0}^{n}\sum_{j=0}^{n}(-2)^{k+j}\binom{n}{k}\binom{n}{j}%
s^{k+1}\int_{0}^{\infty }t^{1+k+j}e^{(iz-1)t}dt \\
&=&\sum_{k=0}^{n}\sum_{j=0}^{n}i^{1+k+j}(-2)^{k+j}\binom{n}{k}\binom{n}{j}%
\Gamma (1+k+j)\frac{s^{k+1}}{(z+i)^{1+k+j}}
\end{eqnarray*}%
Now, since $u_{k,j}(z)=1/(z+i)^{1+k+j}$ is analytic on the upper half plane,
then its maximum on the ball $B(z,\delta )$ is bounded by the average of $%
u_{k,j}(z)$ on the ball. For this reason, $u_{k,j}(z)\in L^{1}(\mathbf{U})$
implies $u_{k,j}(z)\in MB$. As a result, also $\Phi _{n}(z)\in MB$.
\end{proof}

Combining this Lemma with the theorem above, we can assure the existence of
Wavelet frames with windows $\Phi _{n}$.

\begin{theorem}
It is possible to choose a $\delta >0$ such that every uniformly discrete
set $\Gamma $ satisfying $B(z,\delta )\cap \Gamma \neq \varnothing $ for
every $z$, $\mathcal{W}\left( \Phi _{n},\Lambda )\right) $ is a wavelet
frame for $H^{2}(\mathbf{U})$.
\end{theorem}

\begin{corollary}
It is possible to choose a $\delta $ such that every uniformly discrete set $%
\Gamma $ satisfying $B(z,\delta )\cap \Gamma \neq \varnothing $ for every $z$%
, $\Gamma $ is a sampling sequence for $\mathcal{A}^{n}(\mathbf{U})$.
\end{corollary}

\begin{corollary}
It is possible to choose a $\delta $ such that every uniformly discrete set $%
\Gamma $ satisfying $B(z,\delta )\cap \Gamma \neq \varnothing $ for every $z$%
, $\mathcal{W}\left( \Phi _{n},\Lambda )\right) $ is a wavelet superframe
for $H^{2}(\mathbf{U})$.
\end{corollary}

\begin{proof}
The arguments we have used in Lemma 2 can be adapted to prove the existence
of superframes of this form, since
\begin{equation*}
\left\langle \mathbf{\Phi }_{n}\mathbf{,}T_{x}D_{s}\mathbf{\Phi }%
_{n}\right\rangle _{\mathcal{H}}=\sum_{k=0}^{n-1}\left\langle \Phi _{k}%
\mathbf{,}T_{x}D_{s}\Phi _{k}\right\rangle _{H^{2}(\mathbf{U})}\in MB.
\end{equation*}
\end{proof}

\subsection{Necessary conditions on the hyperbolic lattice}

Now we give necessary conditions for the sampling sequences and frames to
exist. Our next result is an upper bound on the size of the parameters $a$
and $b$ (or a lower bound on density) necessary to generate sampling
sequences in the true polyanalytic Bergman space, or wavelet frames with the
functions $\Phi _{n}$. For this purpose, we will adapt the proof of the
necessity part of Theorem 1.1 in \cite{Seip}. An essential step is to
associate to $\Gamma (a,b)$ the analytic function
\begin{equation*}
h(z)=\left( \prod_{k=0}^{\infty }\frac{\sin \pi b^{-1}a^{-k}(ia^{k}-z)}{\sin
\pi b^{-1}a^{-k}(ia^{k}+z)}\right) \left( \prod_{m=1}^{\infty }e^{\frac{2\pi
}{b}}\frac{\sin \pi b^{-1}a^{m}(z-ia^{-m})}{\sin \pi b^{-1}a^{m}(z+ia^{-m})}%
\right) \text{,}
\end{equation*}%
which vanishes in $\Gamma (a,b)$ and plays the role of the sine function in
the Paley-Wiener space and of the Weierstrass $\sigma $-function in the
Bargmann-Fock space. Following \cite{Seip}, one can check that%
\begin{equation*}
h(az)=-e^{-\frac{2\pi }{b}}h(z).
\end{equation*}%
Then, from an estimate of the growth of $h$ on the strip $a^{-1/2}<y<a^{1/2}$%
, the following global estimate is obtained:%
\begin{equation}
\left\vert h(z)\right\vert \leq Cs^{-\frac{2\pi }{b\ln a}}.  \label{growth}
\end{equation}%
Another ingredient in the proof is a result from\textbf{\ \cite{AscensiBruna}%
, }which gives the stability, with respect to the jittered error, for
general wavelet spaces.

\textbf{Theorem C \cite[Theorem 4.4]{AscensiBruna}. }\emph{Let }$W_{g}$\emph{%
\ be the wavelet space associated to the analyzing wavelet }$g$\emph{. If }$%
\Lambda =\{z_{j}\}$\emph{\ is a sampling set for }$W_{g},$\emph{\ there
exists }$\delta >0$\emph{\ such that if }$\Gamma =\{w_{j}\}$\emph{\
satisfies }$\rho (z_{j},w_{j})<\delta $\emph{\ for all }$j$\emph{, then }$%
\Gamma $\emph{\ is also a sampling set.}

\begin{theorem}
If $z_{mk}=a^{m}(bk+i)$ is a sampling sequence for $\mathcal{A}^{n}(\mathbf{U%
})$, then
\begin{equation*}
b\log a<2\pi (n+1).
\end{equation*}
\end{theorem}

\begin{proof}
Suppose that $\Gamma (a,b)$\ is a sampling sequence. The growth estimate (%
\ref{growth}) gives:
\begin{equation*}
h(z)\in A(\mathbf{U})\Longleftrightarrow \frac{4\pi }{b\ln a}<1\text{.}
\end{equation*}%
As a result, if $b\log a>4\pi (n+1)$, then $h^{n+1}(z)\in A(\mathbf{U})$.
Thus, there exists a $f\in L^{2}(%
\mathbb{R}
^{+})$ such that $h^{n+1}(z)=Ber$ $f(z)$. Now consider the function
\begin{equation*}
H(z)=\frac{1}{(2i)^{n}n!}\left( \frac{d}{dz}\right) ^{n}\left[
s^{n}h^{n+1}(z)\right] \text{.}
\end{equation*}%
Clearly,%
\begin{equation*}
H(z)=\frac{1}{(2i)^{n}n!}\left( \frac{d}{dz}\right) ^{n}\left[ s^{n}Berf(z)%
\right] =\emph{Ber}^{n}f(z)\in \mathcal{A}^{n}(\mathbf{U})\text{.}
\end{equation*}%
Since $H(z)$ vanishes on $\Gamma (a,b)$, $\Gamma (a,b)$ cannot be a sampling
sequence for $\mathcal{A}^{n}(\mathbf{U})$. It follows that $b\log a\leq
2\pi (n+1).$

To prove that the inequality is strict, observe that, by Theorem C, there
exists a $\delta >0$ such that if $\Gamma =\{w_{mk}\}$\ satisfies $\rho
(z_{mk},w_{mk})<\delta $ for all $m,k$, then $\Gamma $ is also a sampling
sequence. Thus, if $b\log a=2\pi (n+1)$, we can choose $\delta _{0}$ such
that $w_{mk}=a^{m}(bk+i(1-\delta _{0}))$ satisfies $\rho
(z_{mk},w_{mk})<\delta $ and therefore it is a sampling sequence. This is
impossible by the argument in the previous paragraph, since $%
\{w_{mk}\}=\Gamma (a,b/(1-\delta _{0}))$ and $b/(1-\delta _{0})\log a>b\log
a=2\pi (n+1)$.
\end{proof}

\begin{corollary}
If $\mathcal{W}\left( \Phi _{n},\Lambda )\right) $ is a wavelet frame in $%
H^{2}(\mathbf{U})$, then
\begin{equation*}
b\log a<2\pi (n+1).
\end{equation*}
\end{corollary}

\begin{proof}
The equivalence between the wavelet frame condition and the sampling
condition in the true polyanalytic Bergman space follows from the identity:
\begin{equation*}
\left\langle f\mathbf{,}D\mathbf{_{s}}T_{x}\Phi _{n}\right\rangle _{H^{2}(%
\mathbf{U})}=W_{\Phi _{n}}f(x,s)=s\emph{Ber}^{n}\mathcal{F}f(z)\text{,}
\end{equation*}%
from where it is easily seen that (\ref{frame}) and (\ref{waveframe}) are
equivalent.
\end{proof}

\begin{remark}
Minor adaptations in the proofs of this paper allow us to introduce a weight
in the spaces, by considering the Bergman norm
\begin{equation*}
\int_{\mathbf{U}}\left\vert f(z)\right\vert ^{2}s^{\alpha }dxds<\infty \text{%
,}
\end{equation*}%
with $\alpha >-1$. In this case, $\Lambda $ is a sampling sequence in the
space $\mathcal{A}_{\alpha }^{n}(\mathbf{U})$ if and only if $\mathcal{W}%
(\Phi _{k}^{\alpha },\Lambda )$ is a wavelet frame in $H^{2}(\mathbf{U})$,
where $\Phi _{n}^{\alpha }$ is defined by $\mathcal{F}\Phi _{n}^{\alpha
}(t)=t^{\frac{1}{2}}\emph{l}_{n}^{2\alpha }(2t)$. A necessary condition for
this to happen is:%
\begin{equation*}
b\log a<2\pi \frac{n+1}{\alpha +1}.
\end{equation*}
\end{remark}

\begin{remark}
Since the superframe property requires every system $\mathcal{W}(\Phi
_{k},\Lambda )$ to be a frame, it follows that $b\log a<2\pi $ is a
necessary condition for $\mathcal{W}(\mathbf{\Phi }_{n},\Lambda )$ to be a
wavelet superframe.
\end{remark}

\end{document}